\documentclass[letterpaper,11pt]{article}

\usepackage{ucs}
\usepackage[utf8x]{inputenc}
\usepackage{graphicx}
\usepackage{amsfonts}
\usepackage{dsfont}
\usepackage{amssymb}
\usepackage{amsmath}
\usepackage{amsthm}
\usepackage{enumerate}
\usepackage{stmaryrd}
\usepackage{fullpage}
\usepackage{ifthen}
\usepackage{subfigure}
\usepackage{epic}
\usepackage{authblk}
\usepackage{textcomp}
\usepackage{mathrsfs}
\usepackage{bm}
\usepackage[small]{caption}

\usepackage[hypertexnames=false,colorlinks=true,linkcolor=blue,citecolor=blue]{hyperref}
\usepackage[numbers,comma,square,sort&compress]{natbib}
\usepackage[letterpaper,text={6.5in,9in},centering]{geometry}

\usepackage{color}
\usepackage{titlesec}

\graphicspath{{eps/}{pdf/}}

\newcommand{\bqq}{\begin{equation}}
\newcommand{\eqq}{\end{equation}}
\newcommand{\bqs}{\begin{equation*}}
\newcommand{\eqs}{\end{equation*}}

\newcommand{\C}{\mathbb{C}} 
\newcommand{\I}{\mathbb{I}}
\newcommand{\J}{\mathbb{J}}
\newcommand{\bH}{\mathbb{H}}
\newcommand{\K}{\mathbb{K}}
\newcommand{\R}{\mathbb{R}} 

\newcommand{\W}{\mathbb{W}}

\renewcommand{\O}{\mathcal{O}}

\newcommand{\bG}{\mathbf{G}}
\newcommand{\boldH}{\mathbf{H}}
\newcommand{\cL}{\mathcal{L}}
\newcommand{\vp}{\varphi}
\newcommand{\md}{\mathrm{d}}
\newcommand{\mbi}{\mathbf{i}}

\newtheorem{lem}{Lemma}[section]
\newtheorem{thm}{Theorem}
\newtheorem{prop}[lem]{Proposition}
\newtheorem{cor}[lem]{Corollary}

\newenvironment{Proof}[1][.]%
 {\begin{trivlist}\item[]\textbf{Proof#1 }}%
 {\hspace*{\fill}$\rule{0.3\baselineskip}{0.35\baselineskip}$\end{trivlist}}
 
  {\begin{trivlist}\item[]{\bf Hypothesis #1 }\em}{\end{trivlist}}

\numberwithin{equation}{section}

\title{Asymptotic stability of the critical Fisher-KPP front using pointwise estimates}

\author[1]{Gr\'egory Faye\footnote{email: \texttt{gregory.faye@math.univ-toulouse.fr}}}
\author[2]{Matt Holzer\footnote{email: \texttt{mholzer@gmu.edu}}}
\affil[1]{\small CNRS, UMR 5219, Institut de Math\'ematiques de Toulouse, 31062 Toulouse Cedex, France}
\affil[2]{\small Department of Mathematical Sciences, George Mason University, Fairfax, VA 22030, USA}

\begin{document}
\maketitle

\begin{abstract}
We propose a simple alternative proof of a famous result of Gallay regarding the nonlinear asymptotic stability of the critical front of the Fisher-KPP equation which shows that perturbations of the critical front decay algebraically with rate $t^{-3/2}$ in a weighted $L^\infty$ space. Our proof is based on pointwise semigroup methods and the key remark that the faster algebraic decay rate $t^{-3/2}$ is a consequence of the lack of an embedded zero of the Evans function at the origin for the linearized problem around the critical front.
\end{abstract}



\section{Introduction}

We revisit the asymptotic stability analysis of Gallay \cite{gallay94} for the critical Fisher-KPP front of the following scalar parabolic equation
\bqq
u_t = u_{xx}+cu_x+f(u), \quad t>0, \quad x\in\R,
\label{kpp}
\eqq
where $f:\R\rightarrow\R$ is a $\mathscr{C}^2$ map satisfying $f(0)=f(1)=0$, $f'(0)>0$, $f'(1)<0$ and $f''(u)<0$ for all $u\in(0,1)$ and $c>0$.  For such an example, it is well known that for any wavespeed $c\geq2\sqrt{f'(0)}:=c_*$,  there exist monotone traveling front solutions $q(x)$ connecting $u=1$ at $-\infty$ and $u=0$ at $+\infty$  where the front profile $q$ is solution of the second order ODE
\bqq
0= q_{xx} + c q_x + f(q).
\label{front}
\eqq
The stability of traveling fronts for the Fisher-KPP equation has been studied by many authors.   For the super-critical family of fronts propagating with speeds $c>c^*$, stability was established by Sattinger using exponential weights to stabilize the essential spectrum and yield exponential in time stability; see \cite{sattinger}.   Stability of the critical front was established by \cite{kirchgassner92}, with extensions and refinements achieved in  \cite{bricmont92,eckmann94,gallay94}.  The sharpest of these results for the Fisher-KPP equation is \cite{gallay94}, where perturbations of the critical front are shown to converge in an exponentially weighted $L^\infty$ space with algebraic rate $t^{-3/2}$.  Of course, we also mention that strong results concerning the convergence of compactly supported initial data to traveling fronts are possible for (\ref{kpp}) using comparison principle techniques; see for example \cite{AW78}.

The primary challenge presented by the critical front is that it is not possible to stabilize the essential spectrum using exponential weights.  This is due to the presence of absolute spectrum at $\lambda=0$ in the form of a branch point of the dispersion relation of the asymptotic system near $+\infty$.  The presence of continuous spectrum near the origin suggests algebraic decay and one might further anticipate heat kernel type decay of perturbations.  As we note above, perturbations of the critical front are known to converge slightly faster --  in an exponentially weighted $L^\infty$ space with algebraic rate $t^{-3/2}$; see \cite{gallay94}. 

Our approach is similar to that of \cite{sattinger} where the linear eigenvalue problem is studied in an exponentially weighted space and resolvent estimates are obtained via inverse Laplace transform.  For the super-critical fronts studied in \cite{sattinger} the Laplace inversion contours can be placed in the stable half plane thereby simplifying the analysis.  No such extension is possible here and we instead approach the problem using pointwise semigroup estimates.  Pointwise semigroup methods were introduced by Zumbrun and Howard \cite{zumbrun98} and have been developed over the past several decades to address stability problems where the essential spectrum can not be separated from the imaginary axis.  Applications include stability of viscous shock waves; see \cite{zumbrun98,howard02}, stability and instability of spatially periodic patterns; see \cite{johnson11}, stability of defects in reaction-diffusion equations; see \cite{beck14}, and more recently stability of stationary reaction-diffusion fronts; see \cite{li16}, to mention a few.  

A rough outline of our approach is as follows.
\begin{itemize}
\item Working in an exponentially weighted space with weight $\omega(x)$, we construct bounded solutions $\varphi^\pm(x)$ for the eigenvalue problem $\mathcal{L}p=\lambda p$ on $\mathbb{R}^\pm$ where $\mathcal{L}$ is a linear operator describing the   linearized eigenvalue problem  transformed to the weighted space.
\item Find bounds for the pointwise Green's function, 
\bqs
\bG_\lambda(x,y)=\left\{
\begin{array}{ll}
\dfrac{\vp^+(x)\vp^-(y)}{\W_\lambda(y)}, & x \geq y,\\
& \\
\dfrac{\vp^-(x)\vp^+(y)}{\W_\lambda(y)}, & x \leq y,
\end{array}
\right.
\eqs
where the Wronskian $\W_\lambda(y):=\vp^+(y)\vp^{-'}(y)-\vp^{+'}(y)\vp^{-}(y)$, is often referred to as the Evans function; see \cite{alexander90}.
\item Apply the inverse Laplace transform, and by a suitable choice of inversion contour show that the Green's function 
\bqq
\bG(t,x,y)=\frac{1}{2\pi \mathbf{i}}\int_\Gamma e^{\lambda t} \bG_\lambda(x,y)\mathrm{d}\lambda,
\label{timegreen}
\eqq
decays pointwise with algebraic rate $t^{-3/2}$.  
\item  Apply $L^p$ estimates to the nonlinear solution expressed using  Duhamel's formula,
\bqs
p(t,x) = \int_\R \bG(t,x,y)p_0(y)\mathrm{d}y+\int_0^t \int_\R \bG(t-\tau,x,y) \omega(y) \mathcal{N}(q_*(y),\omega(y)p(\tau,y))p(\tau,y)\mathrm{d}y\mathrm{d}\tau,
\eqs
 to show that the nonlinear system also exhibits the same algebraic decay rate.
\end{itemize}

This approach is motivated by the observation in \cite{sandstede04} that the faster algebraic decay rate is a consequence of the lack of an embedded zero of the Evans function at $\lambda=0$ (the analytic extension of the Evans function to the branch point is possible due to the Gap Lemma; see \cite{gardner98,kapitula98}).  Indeed, the critical front has weak exponential decay near $x=+\infty$, 
\bqq
q_*(x)\underset{+\infty}{\sim}  b x e^{-\gamma_*x},
\label{decayfront}
\eqq
where $\gamma_*:=c_*/2$ and for some $b>0$.  This weak exponential decay implies that the derivative of the wave also has weak exponential decay and therefore does not lead to a zero of $\W_\lambda(y)$ at $\lambda=0$.  In the outline of our argument, this fact comes into play when we require bounds on the supremum of  $\bG(t,x,y)$.  We find that this quantity is  dominated by a region where $\bG_\lambda(x,y)$ resembles $Ce^{-\sqrt{\lambda}(x-y)}$.  Of course, this is exactly the Laplace transform of the derivative of the heat kernel; from which we naturally expect algebraic decay with rate $t^{-3/2}$.  

We now state our main result.  Let $\omega(x)>0$ be a positive, bounded, smooth weight function of the form
\bqs
\omega(x)=\left\{ 
\begin{array}{ll}
e^{-\gamma_*x} & x \geq 1,\\
e^{\beta x} & x \leq -1,
\end{array}
\right.
\eqs
for some $0<\beta<-\frac{c_*}{2}+\sqrt{\frac{c_*^2}{4}-f'(1)}$.  

\begin{thm}\label{thmmain} Consider (\ref{kpp}) with initial data $u(0,x)=q_*(x)+v_0(x)$ satisfying $0\leq u(0,x)\leq 1$.  There exist $C>0$ and $\epsilon > 0$ such that if $ v_0(x)$ satisfies, 
\[  \left|\left| \frac{v_0(\cdot)}{\omega(\cdot)} \right|\right|_{L^\infty} + \left|\left| (1+|\cdot|)\frac{v_0(\cdot)}{\omega(\cdot)} \right|\right|_{L^1}<\epsilon,  \]
then the solution $u(t,x)$ is defined for all time and the  critical front is nonlinearly stable in the sense that
\[ \left|\left| \frac{1}{(1+|\cdot|)} \frac{v(t,\cdot)}{\omega(\cdot)} \right|\right|_{L^\infty}\leq \frac{C \epsilon }{(1+t)^{3/2}},\quad t>0, \]
where $v(t,x):=u(t,x)-q_*(x)$.
\end{thm}

Theorem~\ref{thmmain} recovers the sharp algebraic in time $L^\infty$ decay rate of perturbations of the critical front that was obtained in \cite{gallay94}.   The proof in \cite{gallay94} uses as a weight the derivative of the front profile.  In this weighted space, the linearized operator as $x\to\infty$ is equivalent to the radial Laplacian in three dimensions; for which the fundamental solution possesses algebraic decay rate $t^{-3/2}$.  The nonlinear argument relies on scaling variables and the application of renormalization group techniques.  In comparing Theorem~\ref{thmmain} to the main result in \cite{gallay94} we note small differences in the spatial decay rates of the allowable perturbations and note that the result in \cite{gallay94} is stronger than the one presented here in that the author is able to identify an asymptotic profile for the solution in addition to its decay rate. 

The main novel contribution of our study is to present an alternative proof based upon pointwise semigroup methods and make rigorous  the observation in \cite{sandstede04} that the faster algebraic decay rate is a consequence of the lack of an embedded zero of the Evans function at $\lambda=0$.  We contend that the proof of Theorem~\ref{thmmain} presented here is more elementary than that of \cite{gallay94} as it relies on (rather coarse) ODE estimates, contour integration and a standard nonlinear stability argument avoiding the technical PDE estimates and renormalization group theory of \cite{gallay94}. Furthermore, this alternative method paves the way to tackle a broader class of problems.  For example, one could consider the extended Fisher-KPP equation 
\bqq
u_t = - \gamma u_{xxxx}+u_{xx}+f(u), \quad t>0, \quad x\in\R,
\label{efkpp}
\eqq
where $\gamma>0$ is a small parameter and $f(u)$ is as in \eqref{kpp}. For such an equation, there exists a family of fronts with wavespeed $c\geq c_*(\gamma)$, in the limit $\gamma\rightarrow0$, which were shown to be stable in exponentially weighted spaces \cite{RW01}. It could be possible to adapt the above ideas to prove that the critical front decays algebraically with rate $t^{-3/2}$ in a weighted $L^\infty$ space. The calculations in that case are more involved (a four-dimensional system of ODEs), but the general key ingredients remain unchanged.  Along similar lines, the approach developed in this paper could be used to establish precise stability results for  pulled invasion fronts in systems of reaction-diffusion equations.  For example, refinements of the stability results in \cite{focant98,raugel98} may be achievable.  Our aim in the present paper is to illustrate the key ideas in the simple Fisher-KPP setting \eqref{kpp}, where the analysis is very explicit.

 From the perspective of pointwise semigroup methods, the application here is fairly straightforward.  To reinforce the discussion in the preceding paragraph, we regard this relative simplicity to be a strength of this paper.  One mathematical feature of interest is the presence of the branch point at the origin which prevents the continuation of any contour integrals into the left half of the complex plane.  An important reference in this regard is Howard \cite{howard02} where a marginally stable branch point also arises when considering the stability of degenerate viscous shock waves.  The Fisher-KPP equation being studied here is quite different, but considerable similarities remain between the approach taken here and the one in \cite{howard02}.

The rest of the paper is organized as follows.  In Section~\ref{secprelim}, we set up and study the  linearized eigenvalue problem. In Section~\ref{secptwsGreen}, we derive bounds on the pointwise Green's function $\bG_\lambda(x,y)$.  These estimates are leveraged to obtain bounds on the time Green's function $\bG(t,x,y)$ in Section~\ref{sectimeGreen}.  The nonlinear stability argument is then presented in Section~\ref{secnonlinear}. 

\section{Preliminaries and ODE estimates}\label{secprelim}
In this section, we set up the linear stability problem and begin to construct the pointwise  Green's function $\bG_\lambda(x,y)$.

\begin{figure}[!t]
\centering
\includegraphics[width=1\textwidth]{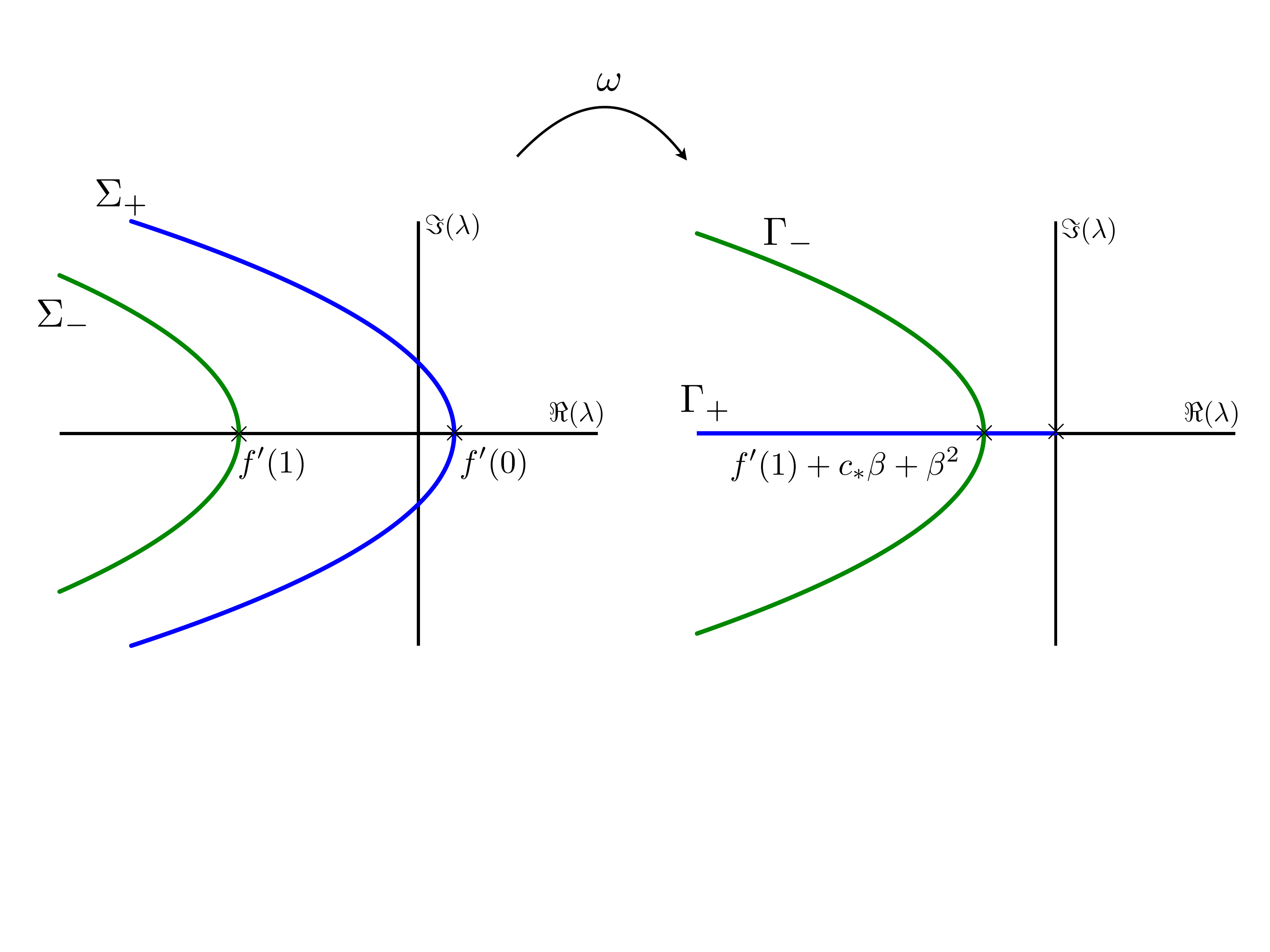}
\caption{Illustration of the action of the weight function $\omega$ on the boundaries of the essential spectrum of the linearized equation around the critical traveling front solution $q_*$. Note that $\Sigma_-$ is mapped to $\Gamma_-$ while $\Sigma_+$ is mapped to the negative real axis $\Gamma_+$ inclusive of the point at $0$.}
\label{fig:ES}
\end{figure}

We work in frame moving to the right with constant velocity $c_*$, so that equation \eqref{kpp} reads
\bqq
u_t = u_{xx}+c_* u_x+f(u), \quad t>0, \quad x\in\R,
\label{kpptw}
\eqq
and writing the solutions $u(t,x)=q_*(x)+v(t,x)$, we obtain the following equation for the perturbation $v(t,x)$:
\bqq
v_t = v_{xx}+c_*v_x+f'(q_*)v+f(q_*+v)-f(q_*)-f'(q_*)v, \quad t>0, \quad x\in\R.
\label{linkpp}
\eqq
Let $\omega(x)>0$ be a positive, bounded, smooth weight function of the form
\bqq
\omega(x)=\left\{ 
\begin{array}{ll}
e^{-\gamma_*x} & x \geq 1,\\
e^{\beta x} & x \leq -1,
\end{array}
\right.
\eqq
for some $\beta>0$ which will be fixed later.  Without loss of generality we assume that $\omega(0)=1$.   We perform a change of variable of the form $v  = \omega  p $, where $p$ now satisfies 
\bqq
p_t = p_{xx}+\left(c_*+2\frac{\omega'}{\omega}\right)p_x+\left(f'(q_*)+c_*\frac{\omega'}{\omega}+\frac{\omega''}{\omega}\right)p+\mathcal{N}(q_*,\omega p)p, \quad t>0, \quad x\in\R,
\label{newkpp}
\eqq
where
\bqs
\mathcal{N}(\mu,\nu):=\frac{1}{\nu}\left(f(\mu+\nu)-f(\mu)-f'(\mu)\nu \right).
\eqs
From now on, we will denote by $\cL$ the linear operator
\bqq
\cL p:= p_{xx}+\left(c_*+2\frac{\omega'}{\omega}\right)p_x+\left(f'(q_*)+c_*\frac{\omega'}{\omega}+\frac{\omega''}{\omega}\right)p,
\label{opL}
\eqq
with dense domain $H^2(\R)$ in $L^2(\R)$. Let us note that for $x\geq1$, the operator $\cL$ reduces to
\bqs
\cL p = p_{xx}+(f'(q_*)-f'(0)) p,
\eqs
while for $x\leq-1$, it becomes
\bqs
\cL p = p_{xx}+(c_*+2\beta)p_x+\left(f'(q_*)+c_*\beta + \beta^2\right)p,
\eqs
from which we impose that $0<\beta<-\frac{c_*}{2}+\sqrt{\frac{c_*^2}{4}-f'(1)}$ so that the essential spectrum of 
\bqs
\cL^-_\infty:=\partial_{xx}+(c_*+2\beta)\partial_x + \left(f'(1)+c_*\beta+\beta^2 \right)
\eqs
 lies to the left of the imaginary axis. We denote by $\Gamma_-$ the parabola in the complex plane defined by the boundary of the essential spectrum of $\cL^-_\infty$, that is 
 \bqs
 \Gamma_-:=\left\{-\ell^2+(c_*+2\beta)\mathbf{i}\ell +f'(1)+c_*\beta+\beta^2 ~|~ \ell \in\R\right\}.
 \eqs
We illustrate in Figure~\ref{fig:ES} the action of the weight function $\omega$ on the boundaries of the essential spectrum of the linearized equation around the critical traveling front solution $q_*$.

For simplicity, we will denote
\bqs
\zeta_1:=c_*+2\frac{\omega'}{\omega}=\left\{
\begin{array}{ll}
0 & x \geq 1,\\
(c_*+2\beta) & x \leq -1,
\end{array}
\right.
 \text{ and } \zeta_0:=f'(q_*)+c_*\frac{\omega'}{\omega}+\frac{\omega''}{\omega}=\left\{
\begin{array}{ll}
f'(q_*)-f'(0) & x \geq 1,\\
f'(q_*)+c_*\beta+\beta^2 & x \leq -1,
\end{array}
\right.
\eqs
with asymptotics
\bqs
\underset{x\rightarrow-\infty}{\lim}\zeta_0(x)=f'(1)+c_*\beta+\beta^2<0 , \quad \text{ and } \quad  \underset{x\rightarrow+\infty}{\lim}\zeta_0(x)=0,
\eqs
such that operator $\cL$ reads
\bqs
\cL p = p_{xx}+\zeta_1p_x+\zeta_0p.
\eqs
The pointwise Green's function  $\bG_\lambda(x,y)$ is a solution of 
\bqq
(\cL-\lambda) \bG_\lambda = -\delta(x-y),
\label{eqgreen}
\eqq
where $\delta(x)$ is the Dirac distribution.  To the right of the essential spectrum solutions can be constructed by identifying exponentially decaying solutions on either half-line and then matching them at $x=y$ enforcing continuity and a jump discontinuity in the derivative,
\[ \lim_{x\to y^-} \frac{\partial \bG_\lambda}{\partial x}+1=\lim_{x\to y^+} \frac{\partial \bG_\lambda}{\partial x}. \]
We begin with the construction of the exponentially decaying solutions of 
\bqq
\cL p = \lambda p.
\label{eigpbm}
\eqq
Our notations will be to denote by $\vp^\pm(x)$ the (unique up to multiplication by a constant) exponentially decaying solutions of \eqref{eigpbm} at $\pm\infty$, and by $\psi^\pm$ a choice of exponentially growing solutions at $\pm\infty$. We can already compute the asymptotic decay and growth rates of $\vp^\pm$ and $\psi^\pm$ from \eqref{eigpbm} at $\pm \infty$. At $+\infty$, we obtain the simpler system
\bqs
p_{xx}=\lambda p,
\eqs
from which we deduce the asymptotic exponential rates $\pm \sqrt{\lambda}$. While at $-\infty$, the system reduces to
\bqs
\cL^-_\infty p = \lambda p,
\eqs
with growth rates given by
\bqs
\mu^\pm(\lambda,\beta):=\frac{-(c_*+2\beta)}{2}\pm\frac{1}{2}\sqrt{(c_*+2\beta)^2-4(f'(1)+c_*\beta+\beta^2)+4\lambda}.
\eqs

Using the above notations, the Green's function $\bG_\lambda(x,y)$ for \eqref{eqgreen} takes the form
\bqs
\bG_\lambda(x,y)=\left\{
\begin{array}{ll}
\dfrac{\vp^+(x)\vp^-(y)}{\W_\lambda(y)}, & x \geq y,\\
& \\
\dfrac{\vp^-(x)\vp^+(y)}{\W_\lambda(y)}, & x \leq y,
\end{array}
\right.
\eqs
where $\W_\lambda(y)$ is the Wronskian,
\bqs
\W_\lambda(y):=\vp^+(y)\vp^{-'}(y)-\vp^{+'}(y)\vp^{-}(y),
\eqs
 and thus satisfies the equation
 \bqs
 \partial_y \W_\lambda(y)=-\zeta_1(y)\W_\lambda(y).
 \eqs
As a consequence and due to the choice of $\omega(0)=1$ we have that
\bqq
 \W_\lambda(y)=\frac{e^{-c_*y}}{\omega^2(y)}\W_\lambda(0), \quad y\in\R.
\label{wronskian}
\eqq
From \eqref{wronskian} and the specific form of the weight function $\omega$, it is then straightforward to check that $\W_\lambda(y)$ simplifies to
\bqs
\W_\lambda(y)=\left\{
\begin{array}{ll}
\W_\lambda(0), & y\geq 1,\\
e^{-(c_*+2\beta)y}\W_\lambda(0), & y\leq -1.
\end{array}
\right.
\eqs

\begin{lem}\label{lem:wronskian}
The Wronskian function $\W_\lambda$ satisfies the following properties:
\begin{itemize}
\item[(i)] $\W_0(y)\neq 0$ for all $y\in\R$;
\item[(ii)] there exists $M_s>0$ such that for all $\lambda$ to the right of $\Gamma_-$ and off the negative real axis with $|\lambda|<M_s$, we have
\bqs
\frac{1}{|\W_\lambda(y)|}\leq C
\eqs
for all $y\in\R$ and some constant $C>0$.
\end{itemize}
\end{lem}

\begin{Proof}
When $\lambda=0$, we have that $\vp^-=\omega^{-1}q_*'$. Using the asymptotic behavior of the critical front at $+\infty$, that is that there exist $a>0$ and $b>0$ such that
\bqs
q_*(y)=(a+by)e^{-\gamma_* y}+\mathcal{O}(y^2e^{-2\gamma_* y}),
\eqs
as $y\rightarrow+\infty$, we deduce that
\bqs
\vp^{-'}(y) \underset{+\infty}{\sim}-\gamma_* b,
\eqs
and as a consequence,
\bqs
\W_0(y)\underset{+\infty}{\sim}-\gamma_* b,
\eqs
which in turn implies that $\W_0(0)\neq0$.  We further note that $\W_\lambda(y)\neq 0 $ for all $\lambda$ to the right of the essential spectrum indicating the absence of unstable point spectrum.  This is a consequence of Theorem 5.5 in \cite{sattinger}.  
We therefore obtain (ii) from the definition of $\W_\lambda$ in \eqref{wronskian}.
\end{Proof}

\begin{lem}\label{lem:phipsi}
Under the assumptions of our main theorem, for the constant $M_s>0$ from Lemma~\ref{lem:wronskian}, and $0<\alpha<\gamma_*$, we have the following estimates on the growth $\psi^\pm$ and decay $\vp^\pm$ modes of \eqref{eigpbm}.
\begin{itemize}
\item[(i)] $(0\leq x)$ For all $|\lambda|\leq M_s$, to the right of $\Gamma_-$ and off the negative real axis, 
\begin{align*}
\vp^+(x)&=e^{-\sqrt{\lambda}x}\left(1+\theta_1^+(x,\lambda)\right),\\
\vp^{+'}(x)&=e^{-\sqrt{\lambda}x}\left(-\sqrt{\lambda}+\theta_2^+(x,\lambda)\right),\\
\psi^+(x)&=e^{\sqrt{\lambda}x}\left(1+\kappa_1^+(x,\lambda)\right),\\
\psi^{+'}(x)&=e^{\sqrt{\lambda}x}\left(\sqrt{\lambda}+\kappa_2^+(x,\lambda)\right),
\end{align*}
where 
\bqs
\theta_1^+(x,\lambda), \kappa_1^+(x,\lambda)=\mathcal{O}(e^{-\alpha x})
\eqs
while 
\bqs
\theta_2^+(x,\lambda), \kappa_2^+(x,\lambda)=\mathcal{O}(\sqrt{\lambda})\mathcal{O}(e^{-\alpha x}).
\eqs
\item[(ii)] $(x\leq 0)$ For all $|\lambda|\leq M_s$ and to the right of $\Gamma_-$, 
\begin{align*}
\vp^-(x)&=e^{\mu^+(\lambda,\beta)x}\left(1+\theta_1^-(x,\lambda)\right),\\
\vp^{-'}(x)&=e^{\mu^+(\lambda,\beta)x}\left(\mu^+(\lambda,\beta)+\theta_2^-(x,\lambda)\right),\\
\psi^-(x)&=e^{\mu^-(\lambda,\beta)x}\left(1+\kappa_1^-(x,\lambda)\right),\\
\psi^{-'}(x)&=e^{\mu^-(\lambda,\beta)x}\left(\mu^-(\lambda,\beta)+\kappa_2^-(x,\lambda)\right),
\end{align*}
where 
\bqs
\theta_{1,2}^-(x,\lambda), \kappa_{1,2}^-(x,\lambda)=\mathcal{O}(e^{\alpha x}).
\eqs
\end{itemize}
The order estimates are uniform in $|\lambda|\leq M_s$.
\end{lem}

\begin{Proof}
As the analysis is similar for each case, we will develop it only for $\vp^+$ and $\psi^+$. Following \cite{zumbrun98}, we first write \eqref{eigpbm} as a first order system of differential equations of the form
\bqq
P'=\mathcal{A}(x,\lambda)P, \quad P=(p,p')^{\mathbf{t}},
\label{eqPsyst}
\eqq
where 
\bqs
\mathcal{A}(x,\lambda):=\left(\begin{matrix} 0 & 1 \\ \lambda-\zeta_0(x) & -\zeta_1(x)  \end{matrix} \right),
\eqs
with associated asymptotic matrices at $x=\pm\infty$
\bqs
\mathcal{A}^+(\lambda):=\left(\begin{matrix} 0 & 1 \\ \lambda & 0  \end{matrix} \right), \quad \text{ and } \quad \mathcal{A}^-(\lambda):=\left(\begin{matrix} 0 & 1 \\ \lambda -(f'(1)+c_*\beta+\beta^2) & -c_*-2\beta  \end{matrix} \right).
\eqs
\paragraph{Case $\vp^+$.} Setting $P(x)=e^{-\sqrt{\lambda}x}Z(x)$, we can rewrite $P'=\mathcal{A}(x,\lambda)P$ as
\bqq
Z' = (\mathcal{A}^+(\lambda)+\sqrt{\lambda}I_2)Z+(\underset{:=\mathcal{B}(x)}{\underbrace{\mathcal{A}(x,\lambda)-\mathcal{A}^+(\lambda)}})Z,
\label{eqZplus}
\eqq
where for $x\geq1$ we have
\bqs
\mathcal{B}(x) = \left(\begin{matrix} 0 & 0 \\ -\zeta_0(x) & 0 \end{matrix} \right), 
\eqs
with $|\zeta_0(x)|=|f'(0)-f'(q_*(x))|=\mathcal{O}(e^{-\alpha x})$ as $x\rightarrow+\infty$. We remark that the matrix $\mathcal{A}^+(\lambda)+\sqrt{\lambda}I_2$ has two eigenvalues: $0$ and $2\sqrt{\lambda}$. As we seek a solution of \eqref{eqZplus} such as $Z(x)\rightarrow Z^+(\lambda)=(1,-\sqrt{\lambda})^{\mathbf{t}}$ as $x\rightarrow+\infty$ we thus have to look for solutions of the integral equation
\bqs
Z(x)=Z^+-\int_x^{+\infty}e^{(\mathcal{A}^+(\lambda)+\sqrt{\lambda}I_2)(x-y)}\mathcal{B}(y)Z(y)\mathrm{d}y,
\eqs
for which we have the explicit formula for $e^{(\mathcal{A}^+(\lambda)+\sqrt{\lambda}I_2)z}$ given by
\bqs
e^{(\mathcal{A}^+(\lambda)+\sqrt{\lambda}I_2)z}=\left(\begin{matrix} \frac{1}{2}+\frac{1}{2}e^{2\sqrt{\lambda}z} & -\frac{1}{2\sqrt{\lambda}}+\frac{1}{2\sqrt{\lambda}}e^{2\sqrt{\lambda}z} \\
-\frac{\sqrt{\lambda}}{2}+\frac{\sqrt{\lambda}}{2}e^{2\sqrt{\lambda}z} & \frac{1}{2}+\frac{1}{2}e^{2\sqrt{\lambda}z}\end{matrix} \right).
\eqs
As a consequence of the specific structure of the matrix $\mathcal{B}(x)$, we obtain a decoupled equation for the first component $z_1(x)$ of $Z(x)$ which we shall obtain as a fixed point of the following map $\mathcal{T}$
\bqs
\mathcal{T}(z_1)(x):=1+\int_x^{+\infty}\left( -\frac{1}{2\sqrt{\lambda}}+\frac{1}{2\sqrt{\lambda}}e^{2\sqrt{\lambda}(x-y)}\right)\zeta_0(y)z_1(y)\md y,
\eqs
for $x\geq1$. The contraction mapping theorem and the remark that $|\zeta_0(x)|=\mathcal{O}(e^{-\alpha x})$ as $x\rightarrow+\infty$ implies the existence of a fixed point of $\mathcal{T}$ on $L^\infty([A,\infty))$ for $A>0$ sufficiently large. Then, we define $\theta_1^+(x,\lambda)$ as
\bqs
\theta_1^+(x,\lambda):=z_1(x)-1.
\eqs
Using an iterative argument, we readily get that $\theta_1^+(x,\lambda)=\mathcal{O}(e^{-\alpha x})$ as $x\rightarrow+\infty$ uniformly in $\lambda$. Furthermore the function $\theta_1^+(\cdot,\lambda)$ is analytic in $\sqrt{\lambda}$, as it can directly be inferred from the fact that
\bqs
\int_x^{+\infty}\left( -\frac{1}{2\sqrt{\lambda}}+\frac{1}{2\sqrt{\lambda}}e^{2\sqrt{\lambda}(x-y)}\right)\zeta_0(y)\md y
 = - \int_x^{+\infty}e^{2\sqrt{\lambda}(x-y)} \left( \int_y^{+\infty}\zeta_0(\tau)\md \tau \right) \md y
 \eqs
 by integration by parts. From the solution $z_1$, we directly obtain an expression for the second component  $z_2$ of $Z(x)$ as
 \bqs
 z_2(x)=-\sqrt{\lambda}+\int_x^{+\infty}\left( \frac{1}{2}+\frac{1}{2}e^{2\sqrt{\lambda}(x-y)}\right)\zeta_0(y)z_1(y)\md y:=-\sqrt{\lambda}+\theta_2^+(x,\lambda),
 \eqs
with the asymptotics $\theta_2^+(x,\lambda)=\mathcal{O}(\sqrt{\lambda})\mathcal{O}(e^{-\alpha x})$ as $x\rightarrow+\infty$. To obtain the conclusion of the lemma and a solution defined for all $x\geq0$, it suffices to flow backward the solution of \eqref{eqZplus} from $x=A$ to $x=0$
. 
\paragraph{Case $\psi^+$.} For $\psi^+$, we perform the change of variable $P(x)=e^{\sqrt{\lambda}x}Z(x)$ in \eqref{eqPsyst}, so that we obtain
\bqq
Z' = (\mathcal{A}^+(\lambda)-\sqrt{\lambda}I_2)Z+\mathcal{B}(x)Z.
\label{eqZpp}
\eqq
We want to construct $\psi^+$ as a solution of \eqref{eqZpp} which satisfies the following two conditions:
\begin{itemize}
\item $\underset{x\rightarrow+\infty}{\lim}Z(x)=\tilde{Z}^+(\lambda)=(1,\sqrt{\lambda})^{\mathbf{t}}$,
\item strong convergence to $\widetilde{Z}^+(\lambda)$.
\end{itemize}
To do so, we first write the variation of constants formula for \eqref{eqZpp}
\bqq
Z(x)=e^{(\mathcal{A}^+(\lambda)-\sqrt{\lambda}I_2)(x-x_0)}Z_0+\int_{x_0}^xe^{(\mathcal{A}^+(\lambda)-\sqrt{\lambda}I_2)(x-y)}\mathcal{B}(y)Z(y)\md y,
\label{eqVC}
\eqq
for some $x_0\geq 1$ and $Z_0\in\C^2$ to be chosen later. The matrix $\mathcal{A}^+(\lambda)-\sqrt{\lambda}I_2$ has two eigenvalues $0$ and $-2\sqrt{\lambda}$ to which we associate a center projection $\Pi_c$ and  a stable projection $\Pi_s$. Noticing that we want to impose convergence of $Z(x)$ as $x\rightarrow+\infty$ along the center direction $(1,\sqrt{\lambda})^{\mathbf{t}}$, we have that
\begin{align*}
Z(x)&=e^{(\mathcal{A}^+(\lambda)-\sqrt{\lambda}I_2)(x-x_0)}\Pi_s Z_0-\int_x^{+\infty} e^{(\mathcal{A}^+(\lambda)-\sqrt{\lambda}I_2)(x-y)}\Pi_c \mathcal{B}(y)Z(y)\md y+ \widetilde{Z}^+(\lambda)\\
&~~~+ \int_{x_0}^xe^{(\mathcal{A}^+(\lambda)-\sqrt{\lambda}I_2)(x-y)}\Pi_s\mathcal{B}(y)Z(y)\md y.
\end{align*}
In order to have strong convergence along the center direction, we further impose
\bqs
e^{-(\mathcal{A}^+(\lambda)-\sqrt{\lambda}I_2)x_0}\Pi_s Z_0+\int_{x_0}^{+\infty}e^{-(\mathcal{A}^+(\lambda)-\sqrt{\lambda}I_2)y}\Pi_s\mathcal{B}(y)Z(y)\md y=0, 
\eqs
such that we look for bounded solutions of
\bqs
Z(x)=\widetilde{Z}^+(\lambda)-\int_x^{+\infty} e^{(\mathcal{A}^+(\lambda)-\sqrt{\lambda}I_2)(x-y)}\mathcal{B}(y)Z(y)\md y.
\eqs
Once again, the first component decouples and we obtain
\bqs
z_1(x) = 1 + \int_x^{+\infty} \left(\frac{1}{2\sqrt{\lambda}}-\frac{1}{2\sqrt{\lambda}}e^{-2\sqrt{\lambda}(x-y)}\right)\zeta_0(y)z_1(y)\md y.
\eqs
We then look for solutions of the form $z_1(x)=1+\kappa_1^+(x,\lambda)$, where $\kappa_1^+(x,\lambda)$ is solution of
\bqs
\kappa(x)=\int_x^{+\infty} \left(\frac{1}{2\sqrt{\lambda}}-\frac{1}{2\sqrt{\lambda}}e^{-2\sqrt{\lambda}(x-y)}\right)\zeta_0(y)\md y+ \int_x^{+\infty} \left(\frac{1}{2\sqrt{\lambda}}-\frac{1}{2\sqrt{\lambda}}e^{-2\sqrt{\lambda}(x-y)}\right)\zeta_0(y)\kappa(y)\md y.
\eqs
By integration by parts and as $|\zeta_0(x)|=\mathcal{O}(e^{-\alpha x})$ as $x\rightarrow+\infty$, we also have that
\bqs
\int_x^{+\infty} \left(\frac{1}{2\sqrt{\lambda}}-\frac{1}{2\sqrt{\lambda}}e^{-2\sqrt{\lambda}(x-y)}\right)\zeta_0(y)\md y=-\int_x^{+\infty}e^{-2\sqrt{\lambda}(x-y)}\left(\int_y^{+\infty}\zeta_0(\tau)\md \tau\right)\md y,
\eqs
for all $\lambda$ such that $2\Re(\sqrt{\lambda})-\alpha<0$. As a consequence, we have the estimate
\bqs
\left| \int_x^{+\infty} \left(\frac{1}{2\sqrt{\lambda}}-\frac{1}{2\sqrt{\lambda}}e^{-2\sqrt{\lambda}(x-y)}\right)\zeta_0(y)\md y \right|=\mathcal{O}(e^{-\alpha x}),
\eqs
which allows us to apply the contraction mapping theorem and get the existence of $\kappa_1^+$ such that $x\mapsto e^{\alpha x} \kappa_1^+(x,\lambda) \in L^\infty([B,\infty))$ for $B>0$ sufficiently large.
\end{Proof}

The following Corollary will be essential in the derivation of our bounds for the pointwise Green's function $\bG_\lambda(x,y)$.  It will turn out that the $\O(e^{-\alpha x})$ bounds for $\kappa_{1,2}^+$ and $\theta_{1,2}^+$ will be insufficient and we will instead require a bound of $\O(\sqrt{\lambda})\O(e^{-\alpha x})$ on their difference.  We remark that a similar cancellation occurs in \cite{howard02} for a different problem where a branch point exists at the origin.  

\begin{cor}\label{cor}
For $\kappa_{1,2}^+$ and $\theta_{1,2}^+$ defined in Lemma \ref{lem:phipsi}, we have that
\bqq
\kappa_{1,2}^+(x,\lambda)-\theta_{1,2}^+(x,\lambda)=\sqrt{\lambda}\Lambda_{1,2}^+(x,\lambda),
\label{eqthka}
\eqq
where $\Lambda_{1,2}^+(x,\lambda)=\mathcal{O}(e^{-\alpha x})$ for $x\geq0$ and $\Lambda_{1,2}^+$ is analytic in $\sqrt{\lambda}$.
\end{cor}

\begin{Proof}
We recall from the previous Lemma that $\theta_1^+$ is solution of
\bqs
\theta(x)=- \int_x^{+\infty}e^{2\sqrt{\lambda}(x-y)} \left( \int_y^{+\infty}\zeta_0(\tau)\md \tau \right) \md y+\int_x^{+\infty}\left( -\frac{1}{2\sqrt{\lambda}}+\frac{1}{2\sqrt{\lambda}}e^{2\sqrt{\lambda}(x-y)}\right)\zeta_0(y)\theta(y)\md y,
\eqs
while $\kappa_1^+$ is solution of
\bqs
\kappa(x)=-\int_x^{+\infty}e^{-2\sqrt{\lambda}(x-y)}\left(\int_y^{+\infty}\zeta_0(\tau)\md \tau\right)\md y+ \int_x^{+\infty} \left(\frac{1}{2\sqrt{\lambda}}-\frac{1}{2\sqrt{\lambda}}e^{-2\sqrt{\lambda}(x-y)}\right)\zeta_0(y)\kappa(y)\md y.
\eqs
As a consequence, we obtain that
\begin{align*}
\kappa_1^+(x,\lambda)-\theta_1^+(x,\lambda)&=-\int_x^{+\infty}\left(e^{-2\sqrt{\lambda}(x-y)}-e^{2\sqrt{\lambda}(x-y)}\right)\left(\int_y^{+\infty}\zeta_0(\tau)\md \tau\right)\md y\\
&~~~+ \int_x^{+\infty} \left(\frac{1}{2\sqrt{\lambda}}-\frac{1}{2\sqrt{\lambda}}e^{-2\sqrt{\lambda}(x-y)}\right)\zeta_0(y)\kappa_1^+(y,\lambda)\md y\\
&~~~- \int_x^{+\infty}\left( -\frac{1}{2\sqrt{\lambda}}+\frac{1}{2\sqrt{\lambda}}e^{2\sqrt{\lambda}(x-y)}\right)\zeta_0(y)\theta_1^+(y,\lambda)\md y.
\end{align*}
Denoting $\Theta_0(y)=-\int_y^{+\infty}\zeta_0(\tau)\md \tau$, we have 
\bqs
\int_x^{+\infty}\left(e^{-2\sqrt{\lambda}(x-y)}-e^{2\sqrt{\lambda}(x-y)}\right)\Theta_0(y)\md y = 2 \sqrt{\lambda} \int_x^{+\infty}\left(e^{-2\sqrt{\lambda}(x-y)}+e^{2\sqrt{\lambda}(x-y)}\right)\left(\int_y^{+\infty}\Theta_0(\tau)\md \tau\right)\md y,
\eqs
where 
\bqs
 \int_x^{+\infty}\left(e^{-2\sqrt{\lambda}(x-y)}+e^{2\sqrt{\lambda}(x-y)}\right)\left(\int_y^{+\infty}\Theta_0(\tau)\md \tau\right)\md y = \mathcal{O}(e^{-\alpha x}),
\eqs
for $x\geq0$. Iterating the argument in the other integral terms, we finally obtain the desired expression \eqref{eqthka}. The proof is then similar for $\kappa_{2}^+(x,\lambda)-\theta_{2}^+(x,\lambda)$.
\end{Proof}

\section{Estimates on the Green's function $\bG_\lambda(x,y)$}\label{secptwsGreen}

Define the following subset of the complex plane,
\bqq \Omega_\delta=\left\{ \lambda\in\mathbb{C}\ \  | \  \mathrm{Re}(\lambda)\geq -\delta_0-\delta_1 |\mathrm{Im}(\lambda)| \right\} \label{eqOmegadelta}, \eqq
where $\delta_{0,1}>0$ are chosen small enough such that $\Gamma_-\cap \Omega_\delta = \emptyset$.

We now derive estimates on $\bG_{\lambda}(x,y)$ in two regimes: sufficiently large $\lambda$ and the remaining values of $\lambda$ near the origin.  

\begin{lem}  \label{lemLbig} There exist some $\eta>0$, $C>0$, and an $M_l>0$ such that if $\lambda\in \Omega_\delta$ and $|\lambda|>M_l$  then  
\bqq \left| \bG_\lambda(x,y)\right| \leq \frac{C}{\sqrt{|\lambda|}} e^{-\sqrt{|\lambda|}\eta  |x-y|}.  \label{eqLbig} \eqq
\end{lem}

\begin{Proof} This is a standard result, see Proposition 7.3 of \cite{zumbrun98}, but we sketch it here for completeness.  To construct $\bG_\lambda$, we require solutions of 
\[ p''+\zeta_1(x)p'+\zeta_0(x)p -\lambda p=0.\]
After scaling the spatial coordinate as $\tilde{x}=\sqrt{|\lambda|}x$, we express this as a first order system,
\bqq
\frac{\md P}{\md \tilde{x}}=\mathcal{A}(\tilde{x},\lambda)P, \quad P=\left(p,q/\sqrt{|\lambda|}\right)^{\mathbf{t}}.
\eqq
Let $\tilde{\lambda}=\frac{\lambda}{|\lambda|}$ at which point we observe
\bqs
\mathcal{A}(\tilde{x},\lambda):=\left(\begin{matrix} 0 & 1 \\ \tilde{\lambda} & 0  \end{matrix} \right)+\frac{1}{\sqrt{|\lambda|}} \left(\begin{matrix} 0 & 0 \\ -\frac{1}{\sqrt{|\lambda|}}\zeta_0\left(\frac{\tilde{x}}{\sqrt{|\lambda|}}\right)& -\zeta_1\left(\frac{\tilde{x}}{\sqrt{|\lambda|}}\right)  \end{matrix} \right).
\eqs
Since $\tilde{\lambda}$ is normalized to lie on the unit circle, there exists an $\eta(\delta)>0$ such that $\Re\left(\sqrt{\tilde{\lambda}} \right)>\eta$ for all $\lambda\in\Omega_\delta$.  Therefore, the leading order system possesses modes that decay at exponential rate $e^{-\sqrt{\tilde{\lambda}}\tilde{x}}$ ($e^{\sqrt{\tilde{\lambda}}\tilde{x}}$) as $\tilde{x}\to\infty$ ($\tilde{x}\to -\infty$).    For $\sqrt{|\lambda|}$ sufficiently large, the full system then has solutions that decay with exponential rate $e^{-\eta \tilde{x}}$ ($e^{\eta \tilde{x}}$).  Reverting to the original variables and imposing continuity and a jump discontinuity in the derivative at $x=y$ we obtain estimates for $\bG_\lambda$ as specified in (\ref{eqLbig}).

\end{Proof}

\begin{lem} \label{lemLsmall}
Under the assumptions of our main theorem and for $|\lambda|\leq M_s$ to the right of $\Gamma_-$ and off the negative real axis, we have the following estimates.
\begin{itemize}
\item[(i)] $y\leq 0 \leq x$
\bqs
\bG_\lambda(x,y)=e^{-\sqrt{\lambda}(x-y)}\mathcal{O}\left(e^{(\mu^+(\lambda,\beta)-\sqrt{\lambda})y} \right);
\eqs
\item[(ii)] $x\leq 0 \leq y$
\bqs
\bG_\lambda(x,y)=e^{\sqrt{\lambda}(x-y)}\mathcal{O}\left( e^{(\mu^+(\lambda,\beta)-\sqrt{\lambda})x} \right);
\eqs
\item[(iii)] $0\leq y \leq x$
\bqs
\bG_\lambda(x,y)=e^{-\sqrt{\lambda}(x-y)}\mathcal{O}\left(e^{(\sqrt{\lambda}-\alpha)y} \right);
\eqs
\item[(iv)] $0\leq x \leq y$
\bqs
\bG_\lambda(x,y)=e^{\sqrt{\lambda}(x-y)}\mathcal{O}\left(e^{(\sqrt{\lambda}-\alpha)x} \right);
\eqs
\item[(v)] $y\leq x \leq 0$
\bqs
\bG_\lambda(x,y)=e^{\mu^-(\lambda,\beta)(x-y)}\mathcal{O}\left(1 \right);
\eqs
\item[(vi)] $x\leq y \leq 0$
\bqs
\bG_\lambda(x,y)=e^{\mu^+(\lambda,\beta)(x-y)}\mathcal{O}\left(1 \right).
\eqs
\end{itemize}
All terms $\mathcal{O}$ are analytic for the $\lambda$ considered here.

\end{lem}

\begin{Proof}
Cases (i)-(ii). For $y\leq 0 \leq x$, we have 
\bqs
\bG_\lambda(x,y)=\frac{\vp^+(x)\vp^-(y)}{\W_\lambda(y)},
\eqs
and according to Lemma \ref{lem:phipsi}, we can write
\bqs
\bG_\lambda(x,y)=\frac{1}{\W_\lambda(y)}e^{-\sqrt{\lambda}x}(1+\theta_1^+(x,\lambda))(1+\theta_1^-(y,\lambda))e^{\mu^+(\lambda,\beta)y}.
\eqs
Now using Lemma \ref{lem:wronskian} and the estimates on $\theta_1^\pm(\cdot,\lambda)$, we conclude that
\bqs
\bG_\lambda(x,y)=e^{-\sqrt{\lambda}(x-y)}\mathcal{O}\left(e^{(\mu^+(\lambda,\beta)-\sqrt{\lambda})y} \right),
\eqs
where the terms containing $\mathcal{O}$ are analytic to the right of $\Gamma_-$ and away from the negative real axis. This concludes the case (i), and (ii) can be treated similarly. 

Cases (iii)-(iv). Let us fix  $0\leq y \leq x$. When $0\leq y$ we need to express $\vp^-(y)$ as a linear combination of $\vp^+(y)$ and $\psi^+(y)$. Thus, we write
\bqs
\vp^-(y)=C(y,\lambda)\vp^+(y)+D(y,\lambda)\psi^+(y).
\eqs
We denote by $\J_\lambda$ and $\I_\lambda$ the following two determinants
\bqs
\J_\lambda(y):=\vp^+(y)\psi^{+'}(y)-\vp^{+'}(y)\psi^{+}(y), \quad \I_\lambda(y):=\vp^-(y)\psi^{+'}(y)-\vp^{-'}(y)\psi^+(y),
\eqs
such that the coefficients $C$ and $D$ are given by
\bqs
C(y,\lambda)=\frac{\I_\lambda(y)}{\J_\lambda(y)},\quad \text{ and } \quad
D(y,\lambda)=\frac{\W_\lambda(y)}{\J_\lambda(y)}.
\eqs
Let us first remark that from Lemma \ref{lem:phipsi}, we have
\bqs
\J_\lambda(y)=(1+\theta_1^+(y,\lambda))(\sqrt{\lambda}+\kappa_2^+(y,\lambda)-(-\sqrt{\lambda}+\theta_2^+(y,\lambda))(1+\kappa_1^+(y,\lambda))=2\sqrt{\lambda}+\mathcal{O}(e^{-\alpha y}),
\eqs
as $y\rightarrow+\infty$. Finally, as both $\vp^+$ and $\psi^+$ are solutions of \eqref{eigpbm}, we have that $\partial_y \J_\lambda(y)=0$ for all $y\geq1$, from which we deduce that
\bqs
\J_\lambda(y)=2\sqrt{\lambda}, \quad \text{ for } y\geq 1.
\eqs
It will be convenient to rewrite $\vp^-(y)$ as
\bqs
\vp^-(y)=\left(C(y,\lambda)+\frac{\W_\lambda(y)}{\J_\lambda(y)e^{-2\sqrt{\lambda}y}}\right)\vp^+(y)+\frac{\W_\lambda(y)}{\J_\lambda(y)}\left(\psi^+(y)-e^{2\sqrt{\lambda}y}\vp^+(y)\right).
\eqs
Using Lemma \ref{lem:phipsi} and the Corollary~\ref{cor}, we see that
\bqs
\psi^+(y)-e^{2\sqrt{\lambda}y}\vp^+(y)=e^{\sqrt{\lambda}y}\left( \kappa_1^+(y,\lambda)-\theta_1^+(y,\lambda)\right)=\sqrt{\lambda}e^{\sqrt{\lambda}y}\Lambda_1^+(y,\lambda),
\eqs
where $\Lambda_1^+(y,\lambda)=\mathcal{O}(e^{-\alpha y})$ for $y\geq0$. Thus, we get
\bqs
\frac{\W_\lambda(y)}{\J_\lambda(y)}\left(\psi^+(y)-e^{2\sqrt{\lambda}y}\vp^+(y)\right)=\W_\lambda(y) \frac{\sqrt{\lambda}}{\J_\lambda(y)}e^{\sqrt{\lambda}y}\Lambda_1^+(y,\lambda)=e^{\sqrt{\lambda}y}\mathcal{O}(e^{-\alpha y}),
\eqs
for all $y\geq1$. On the other hand, we have
\begin{align*}
\frac{\I_\lambda(y)}{\J_\lambda(y)}+\frac{\W_\lambda(y)}{\J_\lambda(y)e^{-2\sqrt{\lambda}y}}&=\frac{1}{\J_\lambda(y)}\left\{ \vp^-(y)\left(\psi^{+'}(y)-e^{2\sqrt{\lambda}y}\vp^{+'} \right) -\vp^{-'}(y)\left(\psi^{+}(y)-e^{2\sqrt{\lambda}y}\vp^{+} \right) \right\}\\
&=\frac{\sqrt{\lambda}}{\J_\lambda(y)}\left\{ \vp^-(y)\Lambda_1^+(y,\lambda) -\vp^{-'}(y)\Lambda_2^+(y,\lambda) \right\}e^{\sqrt{\lambda}y},
\end{align*}
where
\bqs
\vp^-(y)\Lambda_1^+(y,\lambda) -\vp^{-'}(y)\Lambda_2^+(y,\lambda)=\mathcal{O}\left(e^{(\sqrt{\lambda}-\alpha)y}\right).
\eqs
Collecting all these results, we obtain that
\bqs
\bG_\lambda(x,y)=\frac{\vp^+(x)\vp^-(y)}{\W_\lambda(y)}=e^{-\sqrt{\lambda}(x-y)}\mathcal{O}\left(e^{(\sqrt{\lambda}-\alpha)y}\right),
\eqs
for all $0\leq y \leq x$. Case (iv) follows similarly.

Cases (v)-(vi). Let us fix  $y \leq x \leq 0$. In that case, we need to express $\vp^+(x)$ as a linear combination of $\vp^-(x)$ and $\psi^-(x)$. Thus we write
\bqs
\vp^+(x)=A(x,\lambda)\vp^-(x)+B(x,\lambda)\psi^-(x).
\eqs
We denote by $\bH_\lambda$ and $\K_\lambda$ the following two determinants 
\bqs
\bH_\lambda(x):=\vp^-(x)\psi^{-'}(x)-\vp^{-'}(x)\psi^{-}(x), \quad \K_\lambda(x):=\vp^+(x)\psi^{-'}(x)-\vp^{+'}(x)\psi^-(x),
\eqs
such that the coefficients $A$ and $B$ are given by
\bqs
A(x,\lambda)=\frac{\K_\lambda(x)}{\bH_\lambda(x)},\quad \text{ and } \quad
B(x,\lambda)=-\frac{\W_\lambda(x)}{\bH_\lambda(x)}.
\eqs
Once again, using Lemma \ref{lem:phipsi}, we obtain
\bqs
\bH_\lambda(x)=e^{-(c_*+2\beta)x}\left((1+\theta_1^-(x,\lambda))(\mu^-(\lambda,\beta)+\kappa_2^-(x,\lambda))-(\mu^+(\lambda,\beta)+\theta_2^-(x,\lambda))(1+\kappa_1^-(x,\lambda)) \right),
\eqs
from which we deduce that
\bqs
\bH_\lambda(x)=e^{-(c_*+2\beta)x}\mathcal{O}(1),
\eqs
for all $x\geq 0$. Recalling that the Wronskian $\W_\lambda$ can be simplified for $x\leq-1$ to
\bqs
\W_\lambda(x)=e^{-(c_*+2\beta)x}\W_\lambda(0)
\eqs
and owing to Lemma \ref{lem:wronskian}, we have
\bqs
B(x,\lambda)=-\frac{\W_\lambda(x)}{\bH_\lambda(x)}=\mathcal{O}(1)
\eqs
for all $x\leq 0$. On the other hand, we have that
\bqs
\partial_x \K_\lambda(x) = -\zeta_1(x) \K_\lambda(x),
\eqs
such that 
\bqs
\K_\lambda (x) = e^{-\int_0^x\zeta_1(\tau)\md \tau}\K_\lambda(0).
\eqs
From which, we deduce that for all $x\leq-1$ we further have
\bqs
\K_\lambda (x) = e^{-(c_*+2\beta)x}\K_\lambda(0),
\eqs
such that we deduce
\bqs
A(x,\lambda)=\frac{\K_\lambda(x)}{\bH_\lambda(x)}=\mathcal{O}(1),
\eqs
which holds true for all $x\leq 0$. Finally, we get
\bqs
\frac{A(x,\lambda)}{\W_\lambda(y)}\vp^-(x)\vp^-(y)= \frac{e^{\mu^+(\lambda,\beta)(x+y)}}{e^{(\mu^+(\lambda,\beta)+\mu^-(\lambda,\beta))y}}\mathcal{O}(1)=e^{\mu^-(\lambda,\beta)(x-y)}\mathcal{O}\left(e^{(\mu^+(\lambda,\beta)-\mu^-(\lambda,\beta))x} \right),
\eqs
and
\bqs
\frac{B(x,\lambda)}{\W_\lambda(y)}\psi^-(x)\vp^-(y)= \frac{\mu^-(\lambda,\beta)x+e^{\mu^+(\lambda,\beta)y}}{e^{(\mu^+(\lambda,\beta)+\mu^-(\lambda,\beta))y}}\mathcal{O}(1)=e^{\mu^-(\lambda,\beta)(x-y)}\mathcal{O}\left(1\right).
\eqs
Case (vi) follows along similar lines and the proof of the lemma is thereby complete.
\end{Proof}

\section{Estimates on the Green's function $\bG(t,x,y)$}\label{sectimeGreen}

We now use the estimates on the pointwise Green's function $\bG_\lambda(x,y)$ to derive bounds on the temporal Green's function $\bG(t,x,y)$.  

\begin{prop}\label{timeGreen}
Under the assumptions of our main theorem, and for some constants $\kappa>0$, $r>0$ and $C>0$, the Green's function $\bG(t,x,y)$ for $\partial_t p = \cL p$ satisfies the following estimates.
\begin{itemize}
\item[(i)] For $|x-y|\geq Kt$ or $t<1$,  with $K$ sufficiently large, 
\bqs
|\bG(t,x,y)|\leq C\frac{1}{t^{1/2}}e^{-\frac{|x-y|^2}{\kappa t}}.
\eqs
\item[(ii)] For $|x-y|\leq Kt$ and $t\geq 1$, with $K$ as above, 
\bqs
|\bG(t,x,y)|\leq C\left(\frac{1+|x-y|}{t^{3/2}} \right)e^{-\frac{|x-y|^2}{\kappa t}}+Ce^{-rt}.
\eqs
\end{itemize}
\end{prop}
\begin{Proof} 
Case (i): The proof of this case follows as in \cite{zumbrun98}.  We select a contour of integration consisting of the concatenation of a parabolic contour contained in the region where the estimates of Lemma~\ref{lemLbig} are valid and a ray extending to infinity.  Consider the region $\Omega_\delta$ with $\delta_0=0$ for simplicity.  We will slightly abuse notation by setting $\delta_1=\delta>0$ into the definition \eqref{eqOmegadelta} of $\Omega_\delta$. Define $\Gamma_\delta=\left\{-\delta |\ell| + \mathbf{i} \ell ~|~ \ell \in\R \right\}$ as the boundary of this region.  We will establish the necessary estimate by appropriate choice of contours on which to apply the inverse Laplace transform.  Define a secondary contour $\Gamma_\rho$ defined for those values of $\lambda$ for which 
\[ \sqrt{\lambda}=\rho+\mbi k \ , \ k\in\mathbb{R},\] 
implying
\[ \lambda(k)=\rho^2-k^2+2\mbi \rho k. \]
During the course of our proof, $\rho$ is chosen sufficiently large so that when $\lambda(k)\in  \Omega_\delta$ estimate (\ref{eqLbig}) applies.  Define the contours $\Delta_1=\Gamma_\rho\cap \Omega_\delta$ and $\Delta_2$ to be the continuation of the contour to $\infty$ along the contour $\Gamma_\delta$. Then,
\bqs
\left|\bG(t,x,y)\right|\leq\frac{1}{2\pi}\left|\int_{\Delta_1} e^{\lambda t} \bG_\lambda(x,y)\mathrm{d}\lambda\right|+\frac{1}{2\pi}\left|\int_{\Delta_2} e^{\lambda t} \bG_\lambda(x,y)\mathrm{d}\lambda\right|.
\eqs
We apply the estimate from Lemma~\ref{lemLbig} to the integral along the contour $\Delta_1$.  Here we have 
\bqs \frac{1}{2\pi}\left|\int_{\Delta_1} e^{\lambda t} \bG_\lambda(x,y)\mathrm{d}\lambda\right|\leq C  \left|\int_{\Delta_1} e^{\rho^2 t-k^2t}e^{-\sqrt{\rho^2+k^2}\eta |x-y|} \mathrm{d}k\right|,
\eqs
where we have used that $\frac{|\sqrt{\lambda}|}{\sqrt{|\lambda|}}=1$.  Using the bound, $\sqrt{\rho^2+k^2}\geq \rho$ we then factor 
\bqs \frac{1}{2\pi}\left|\int_{\Delta_1} e^{\lambda t} \bG_\lambda(x,y)\mathrm{d}\lambda\right|\leq C e^{\rho^2 t-\rho\eta |x-y|} \left|\int_{\Delta_1} e^{-k^2t} \mathrm{d}k\right|.
\eqs
Now, we must select $\rho$ so that 
\[ \rho^2 t-\rho\eta |x-y|=-\frac{|x-y|^2}{\kappa t},\]
for some $\kappa>0$.  Solving for $\rho$, we obtain
\[ \rho=\frac{|x-y|}{t} \left( \frac{\eta}{2}\pm \frac{1}{2}\sqrt{\eta^2-\frac{4}{\kappa}}\right). \]
Recall that $\eta$ is fixed by Lemma~\ref{lemLbig}, so we may select $\kappa$ sufficiently large  so that the equation for $\rho$ has real roots.  With $\kappa$ fixed, we may now restrict to  $\frac{|x-y|}{t}$ sufficiently large so that the contour $\Delta_1$ lies within the region where the estimate (\ref{eqLbig}) applies and we obtain the bound 
\bqs \frac{1}{2\pi}\left|\int_{\Delta_1} e^{\lambda t} \bG_\lambda(x,y)\mathrm{d}\lambda\right|\leq \frac{C}{\sqrt{t}}  e^{-\frac{|x-y|^2}{\kappa t}}. \eqs
We now turn our attention to the contribution from the integral along $\Delta_2$, focusing on the portion in the upper half plane; the case in the lower half plane is similar.  In this region, the contour can be described as 
\[ \lambda(\ell)=-\delta \ell+\mbi \ell.\]
Intersection points of $\Delta_1$ and $\Delta_2$  satisfy
\[ \rho^2-k^2=-\delta \ell,  \quad 2\rho k=\ell,\]
from which we find
\bqq \label{eqkstar} k^*=\rho\left(\delta+\sqrt{\delta^2+\rho^2}\right), \ \ell_*=2\rho^2\left(\delta+\sqrt{\delta^2+\rho^2}\right). \eqq
We then have 
\bqq \frac{1}{2\pi}\left|\int_{\Delta_2} e^{\lambda t} \bG_\lambda(x,y)\mathrm{d}\lambda\right|\leq C \int_{\ell_*}^\infty\frac{1}{\sqrt{\ell}} e^{-\delta \ell t}d\ell \leq C \int_{4\delta\rho^2}^\infty\frac{1}{\sqrt{\ell}} e^{-\delta \ell t}d\ell \leq \frac{C}{\sqrt{t}} \int_{2\delta\rho \sqrt{t}}^\infty e^{-z^2}dz \leq \frac{C}{\sqrt{t}} e^{-(2\delta\rho \sqrt{t})^2}. \label{eqGtbigL}\eqq
Substituting for $\rho$ and potentially taking a larger  value for $\kappa$, we obtain the bound
\bqs \frac{1}{2\pi}\left|\int_{\Delta_2} e^{\lambda t} \bG_\lambda(x,y)\mathrm{d}\lambda\right|\leq \frac{C}{\sqrt{t}}  e^{-\frac{|x-y|^2}{\kappa t}}, \eqs
and the result follows in this case.  

Finally, parabolic regularity gives the same estimate for short time $t<1$.

Case (ii): In the regime $|x-y|\leq Kt$ and $t\geq1$, the analysis must be carried out in pieces according to the relative positions of $x$ and $y$.  The cases where one, or both, of $x$ and $y$ are positive are similar so we consider the case of $y\leq 0\leq x$ only. Here, we will consider the region $\Omega_\delta$ with $\delta_{0,1}>0$ and set $-\delta_1=\cot(\theta)$ for some some fixed $\pi/2<\theta<\pi$, so that it's boundary $\Gamma_\delta$ can be parametrized by $\Gamma_\delta=\left\{ \lambda(\ell)=-\delta_0+\cos(\theta)|\ell|+\mbi \sin(\theta)\ell ~|~ \ell \in \R \right\}$.

We first summarize our approach.  We will deform the Laplace inversion contour into two pieces -- a parabolic segment surrounding the branch point at $\lambda=0$ and a ray tending to infinity along $\Gamma_\delta$.  Recall the estimate in  Lemma~\ref{lemLsmall},
\[ \bG_\lambda(x,y)=e^{-\sqrt{\lambda}(x-y)}\mathcal{O}\left(e^{(\mu^+(\lambda,\beta)-\sqrt{\lambda})y} \right). \]
We select the parabolic contour in such a way that it is contained in the region where 
\bqq \Re\left( \mu^+(\lambda,\beta)-\sqrt{\lambda}\right)>0. \label{eqremcond} \eqq
  The contour integration can then be divided into four integrals depending on the relevant estimates in each region, see Figure~\ref{fig:contour} for an illustration,
\begin{itemize}
\item $\Gamma_1$ -- a parabolic contour near $\lambda=0$ for which  condition (\ref{eqremcond}) holds  ,
\item $\Gamma_2$ -- a continuation of the contour along the ray $\Gamma_\delta$ for those $\lambda$ values where condition (\ref{eqremcond}) holds  ,
\item $\Gamma_3$ -- a continuation of the contour along the ray $\Gamma_\delta$ for those $\lambda$ values where condition (\ref{eqremcond}) fails, but for which the large $\lambda$ estimates from Lemma~\ref{lemLbig} do not apply  ,
\item $\Gamma_4$ -- a continuation of the contour along the ray $\Gamma_\delta$ for those $\lambda$ values where the large $\lambda$ estimates from Lemma~\ref{lemLbig}  apply.
\end{itemize}

\begin{figure}[!t]
\centering
\includegraphics[width=0.85\textwidth]{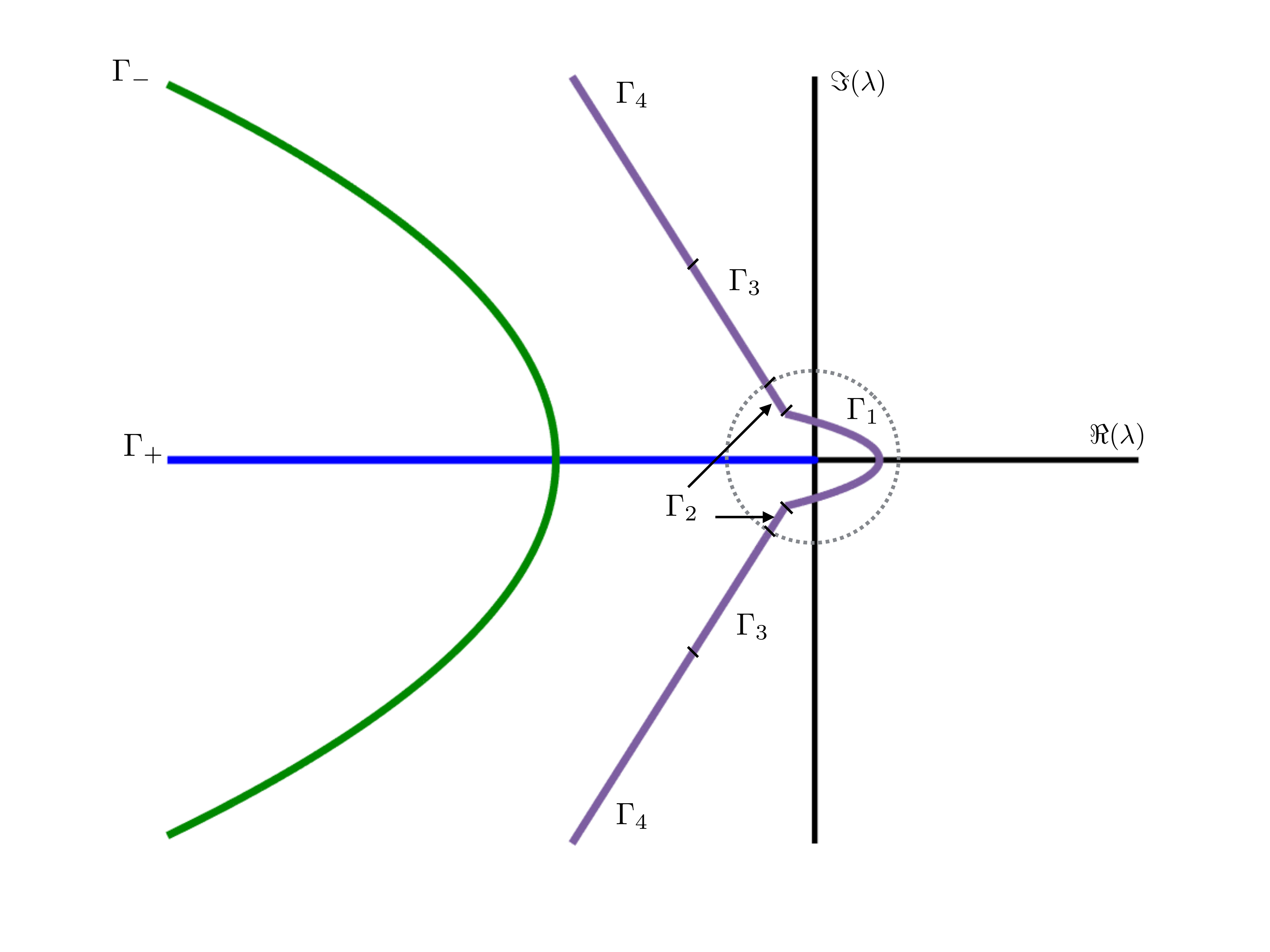}
\caption{Visualization of the 4 contours of integration $\Gamma_i$, $i=1,\ldots,4$, in the regime $|x-y|\leq K t$. $\Gamma_1$ is a parabolic contour near $\lambda=0$ for which  condition (\ref{eqremcond}) holds. $\Gamma_2$ is a continuation of the contour along the ray $\Gamma_\delta$ for those $\lambda$ values where condition (\ref{eqremcond}) holds. $\Gamma_3$ is a continuation of the contour along the ray $\Gamma_\delta$ for those $\lambda$ values where condition (\ref{eqremcond}) fails, but for which the large $\lambda$ estimates from Lemma~\ref{lemLbig} do not apply where we simply used the uniform boundedness of the Green's function. $\Gamma_4$ is a continuation of the contour along the ray $\Gamma_\delta$ for those $\lambda$ values where the large $\lambda$ estimates from Lemma~\ref{lemLbig}  apply.}
\label{fig:contour}
\end{figure}

Note that at $\lambda=0$ we have that $\mu^+(\lambda,\beta)>\sqrt{\lambda}$.  Then there exists $M_\mu>0$ such that for all $\lambda\in\Omega_\delta$, off the imaginary axis and with $|\lambda|<M_\mu$ we have that $\mu^+(\lambda,\beta)>\sqrt{\lambda}$.  

Define the contour $\Gamma_1$ via 
\bqs
\sqrt{\lambda} = \rho+\mbi k,
\eqs
for some $\rho>0$ and for all $k$ such that $\sqrt{\lambda}\in\Omega_\delta$.  

For $\Gamma_2$ through $\Gamma_4$ we work with $\lambda$ in the upper half plane without loss of generality -- estimates for the lower half plane are analogous.  The contours $\Gamma_2$ through $\Gamma_4$ are all given by  
\[ \lambda=-\delta_0 +\cos(\theta) \ell+\mbi \sin(\theta)\ell, \]
for some fixed $\pi/2<\theta<\pi$ and varying intervals of $\ell$.   For any $\rho$ we can compute the intersection points
\bqq k^*(\rho)=-\rho \cot(\theta)+\sqrt{\rho^2 \csc^2(\theta)+\delta_0},\ \quad \ell_1(\rho)=\frac{2\rho k^*(\rho)}{\sin(\theta)}.\label{eql1} \eqq

Recall the estimates in Lemma~\ref{lemLbig} and Lemma~\ref{lemLsmall}.  Write 
\[ \bG_\lambda(x,y)=e^{-\sqrt{\lambda}(x-y)}e^{(\mu^+(\lambda,\beta)-\sqrt{\lambda})x} \boldH_\lambda(x,y),\]
where $\boldH_\lambda(x,y)$ is analytic as a function of $\sqrt{\lambda}$ for $\lambda\in\Omega_\delta$ and uniformly bounded in $(x,y)$.  

For the parabolic contour, we recall that we are interested in $|x-y|\leq Kt$ so we select 
\[ \rho=\frac{|x-y|}{Lt}\]
for some $L$ sufficiently large so that the parabolic contour  $\Gamma_1$ lies within the region for which $\mu^+(\lambda,\beta)>\sqrt{\lambda}$.  Let us note that
\bqs
\md \lambda = 2\mbi \left(\frac{(x-y)}{Lt}+\mbi k\right)\md k \text{ and } \lambda = \frac{(x-y)^2}{L^2t^2}-k^2+2\frac{(x-y)}{Lt}\mbi k.
\eqs
Then
\begin{align} \frac{1}{2\pi \mbi }\int_{\Gamma_1} e^{\lambda t} \bG_\lambda(x,y)\mathrm{d}\lambda  &=  
 \frac{1}{\pi  }e^{\rho^2 t-\rho(x-y)} \int_{-k^*}^{k^*} e^{-k^2 t}e^{2\mbi \rho k -\mbi k(x-y)} \boldH_{\lambda(k)}(x,y)(\rho+\mbi k) \mathrm{d}k \nonumber  \\
&=\frac{1}{\pi  }e^{\rho^2 t-\rho(x-y)} \int_{-k^*}^{k^*} e^{-k^2 t} \left(H_R(x,y, k)+\mbi H_I(x,y,k)\right)(\rho+\mbi k) \mathrm{d}k.
 \end{align}
Here we have expanded $e^{2\mbi \rho k -\mbi k(x-y)} \boldH_{\lambda(k)}$ into its real and imaginary parts.  Since $\bG_\lambda$ is holomorphic and the contour is symmetric with respect to the real axis, we have that $H_R(x,y,k)$ is even in $k$ while $H_I(x,y,k)$ is odd.  Using this, the integral reduces to
\[ \frac{1}{2\pi \mbi }\int_{\Gamma_1} e^{\lambda t} \bG_\lambda(x,y)\mathrm{d}\lambda=\frac{1}{\pi  }e^{\rho^2 t-\rho(x-y)} \int_{-k^*}^{k^*} e^{-k^2 t} \left(\rho H_R(x,y, k)-k H_I(x,y,k)\right) \mathrm{d}k. \]
Boundedness of $H_R$ implies that 
\[ \left|\frac{1}{\pi  }e^{\rho^2 t-\rho(x-y)} \int_{-k^*}^{k^*} e^{-k^2 t} \rho H_R(x,y, k) \mathrm{d}k\right|\leq C\frac{\rho}{\sqrt{t}} e^{\rho^2 t-\rho(x-y)}.   \]
On the other hand, since $H_I(x,y,k)$ is odd in $k$, it can be expressed as $H_I(x,y,k)=k\tilde{H}_I(x,y,k)$, where $\tilde{H}_I(x,y,k)$ is again bounded.  Therefore,  
\[ \left|\frac{1}{\pi  }e^{\rho^2 t-\rho(x-y)} \int_{-k^*}^{k^*} e^{-k^2 t} k^2 \tilde{H}_I(x,y, k) \mathrm{d}k\right| \leq \frac{C}{t^{3/2}} e^{\rho^2 t-\rho(x-y)}.   \]
Using our choice of $\rho$, we then arrive at the estimate 
\[  \frac{1}{2\pi }\left|\int_{\Gamma_1} e^{\lambda t} \bG_\lambda(x,y)\mathrm{d}\lambda\right| \leq C\frac{1+|x-y|}{t^{3/2}} e^{-\frac{|x-y|^2}{\kappa t}},\]
for $\kappa=L^2/(L-1)$.

Next consider the integral along $\Gamma_2$.  Here we have
\begin{align} \frac{1}{2\pi}\left|\int_{\Gamma_2} e^{\lambda t} \bG_\lambda(x,y)\mathrm{d}\lambda\right|  &\leq   
 Ce^{-\delta_0 t} \int_{\ell_1}^{\ell_2} e^{t\cos(\theta)\ell -(x-y)\mathrm{Re}(\sqrt{-\delta_0+\cos(\theta)\ell+\mbi \sin(\theta)\ell}) }\mathrm{d}\ell \nonumber \\
&\leq Ce^{-\delta_0 t} \int_{\ell_1}^{\ell_2} e^{t\cos(\theta)\ell } \mathrm{d}\ell \leq C \frac{e^{-\delta_0 t}}{t}e^{t\cos(\theta)\ell_1}  ,
\end{align}
we were have used that $(x-y)\mathrm{Re}(\sqrt{-\delta_0+\cos(\theta)\ell+\mbi \sin(\theta)\ell}) >0$ and subsequently integrated noting that $\cos(\theta)<0$.    Now, by virtue of (\ref{eql1}) we observe that 
\[ t\cos(\theta) \ell_1(\rho) < 2t\cot(\theta)\rho^2\left( -\cot(\theta)+\csc(\theta)\right)<0.\]
We therefore obtain the bound
\[ \frac{1}{2\pi}\left|\int_{\Gamma_2} e^{\lambda t} \bG_\lambda(x,y)\mathrm{d}\lambda\right|  \leq C\frac{e^{-\delta_0 t}}{t} e^{-\frac{|x-y|^2}{\kappa t}} \leq C\frac{1}{t^{3/2}} e^{-\frac{|x-y|^2}{\kappa t}},\]
for some $\kappa>0$.

Next consider the integral along $\Gamma_3$.  The key distinction in this case is that (\ref{eqremcond}) does not hold.  We instead  use the fact that $\bG_\lambda(x,y)$ is uniformly bounded in this region to obtain
\[ \frac{1}{2\pi }\left|\int_{\Gamma_3} e^{\lambda t} \bG_\lambda(x,y)\mathrm{d}\lambda\right| \leq Ce^{-rt}, \]  
for some $r>0$.  Finally, we consider the integral along $\Gamma_4$.   The large $\lambda$ bounds apply here and the analysis follows as in (\ref{eqGtbigL}).  We find 
\[ \frac{1}{2\pi }\left|\int_{\Gamma_4} e^{\lambda t} \bG_\lambda(x,y)\mathrm{d}\lambda\right| \leq C e^{-\delta_0 t} \int_{\ell_3}^\infty \frac{1}{\sqrt{\ell}} e^{\cos(\theta)\ell t}\mathrm{d}\ell\leq C e^{-rt} , \]  
for some $r>0$.

This concludes the proof of case (ii) with $y\leq 0\leq x$.  To complete the proof similar estimates are required to the five other orderings as in Lemma~\ref{lemLsmall}.  Since the proofs are analogous to the present case, we omit the details.  Note that for cases (v) and (vi) in Lemma~\ref{lemLsmall} we could achieve sharper bounds, but these are not required for the proof of Theorem~\ref{thmmain} so we do not pursue this here. 

\end{Proof}

In the following section, we will require bounds on integrals of the form,
\[ \left| \int_\R \bG(t,x,y)h(y)\mathrm{d}y \right|, \]
for some $h(y)$.

\textbf{Short time estimate.} We apply Proposition~\ref{timeGreen}(i) for $t < 1$ to obtain the usual short time estimate,
\bqq 
\left| \int_\R \bG(t,x,y)h(y)\mathrm{d}y \right| \leq C \| h\|_{L^\infty}.
\label{eqGestsmall}
\eqq

\textbf{Large time estimate.} Applying Proposition~\ref{timeGreen} for $t\geq 1$, we can then estimate 
\begin{align}
\left| \int_\R \bG(t,x,y)h(y)\mathrm{d}y \right| &\leq \int_{-\infty}^{x-Kt} C\frac{1}{t^{1/2}}e^{-\frac{|x-y|^2}{\kappa t}}|h(y)|\mathrm{d}y \nonumber \\
&+\int_{x-Kt}^{x+Kt} C\left( \left(\frac{1+|x-y|}{t^{3/2}} \right)e^{-\frac{|x-y|^2}{\kappa t}}+e^{-rt}\right)|h(y)|\mathrm{d}y  \nonumber \\
&+\int_{x+Kt}^{\infty }C\frac{1}{t^{1/2}}e^{-\frac{|x-y|^2}{\kappa t}}|h(y)|\mathrm{d}y.
\end{align}
Note that the terms $\frac{1}{t^{1/2}}e^{-\frac{|x-y|^2}{\kappa t}}$ in the first and last integrals are maximized at the boundary and therefore the contribution from these integrals decay exponentially in time.  For the middle integral we take the limits of integration to infinity and use $1+|x-y|\leq 1+|x|+|y|+|x||y|$ from which we obtain
\[ \left| \int_\R \bG(t,x,y)h(y)\mathrm{d}y \right| \leq Ce^{-rt}\int_{\mathbb{R}}(1+|y|)|h(y)|\mathrm{d}y +C \frac{1+|x|}{(1+t)^{3/2}}\int_\R(1+|y|)|h(y)|\md y.\]
Due to the exponential decay of the first integral, we in fact have 
\bqq \left| \int_\R \bG(t,x,y)h(y)\mathrm{d}y \right| \leq C \frac{1+|x|}{(1+t)^{3/2}}\int_\R(1+|y|)|h(y)|\md y.\label{eqGest} \eqq

\section{Nonlinear Stability: Proof of Theorem~\ref{thmmain}}\label{secnonlinear}

We recall that we look for perturbations $p(t,x)$ which satisfy the evolution equation 
\bqq
p_t = p_{xx}+\left(c_*+2\frac{\omega'}{\omega}\right)p_x+\left(f'(q_*)+c_*\frac{\omega'}{\omega}+\frac{\omega''}{\omega}\right)p+\mathcal{N}(q_*,\omega p)p, \quad t>0, \quad x\in\R.
\label{evolp}
\eqq
The Cauchy problem associated to equation \eqref{evolp} with initial condition $p_0 \in L^1(\R)\cap L^\infty(\R)$ with $\int_\R |y| |p_0(y)|\md y < +\infty$ is locally well-posed in $L^\infty(\R)$. As a consequence, we let $T_*>0$ be the maximal time of existence of a solution $p\in L^\infty(\R)$ with initial condition $p_0 \in L^1(\R)\cap L^\infty(\R)$ with $\int_\R |y| |p_0(y)|\md y < +\infty$ of the associated integral formulation of \eqref{evolp} given by
\bqq
p(t,x) = \int_\R \bG(t,x,y)p_0(y)\mathrm{d}y+\int_0^t \int_\R \bG(t-\tau,x,y) \mathcal{N}(q_*(y),\omega(y)p(\tau,y))p(\tau,y)\mathrm{d}y\mathrm{d}\tau.
\label{inteqp}
\eqq
For $t\in[0,T_*)$, we define
\bqs
\Theta(t) = \underset{0\leq \tau \leq t}{\sup}~\underset{x\in\R}{\sup} ~ (1+\tau )^{3/2} \frac{|p(\tau,x)|}{1+|x|}.
\eqs
Furthermore, by our regularity assumption on $f$, there exists a positive nondecreasing function $\chi:\R_+\rightarrow\R_+$, such that
\bqs
| \mathcal{N}(q_*,\omega p )p | \leq \chi (R) \omega p^2, \quad |\omega p| \leq R.
\eqs 
Now, using the fact that $\omega$ is exponentially localized, we have that for all $0\leq \tau \leq t$ and all $y\in\R$
\bqs
|\mathcal{N}(q_*(y),\omega(y)p(\tau,y))p(\tau,y)| \leq \chi( \omega_\infty \Theta(t) ) \omega(y) p^2(\tau,y) \leq \frac{1}{(1+\tau)^3}\chi( \omega_\infty \Theta(t) ) (1+|y|)^2 \omega(y) \Theta(t)^2,
\eqs
where $\omega_\infty:=\underset{x\in\R}{\sup}\left((1+|x|)\omega(x) \right)<\infty$.

We now bound each term of the right hand side of \eqref{inteqp} for $0\leq t<1$ and $t\geq 1$. If $0<t<1$, we use the estimate \eqref{eqGestsmall} to obtain
\bqs \left| \int_\R \bG(t,x,y)p_0(y)\mathrm{d}y \right| \leq C \|p_0\|_{L^\infty}\eqs 
together with
\bqs
\left| \int_\R \bG(t-\tau,x,y) \mathcal{N}(q_*(y),\omega(y)p(\tau,y))p(\tau,y)\mathrm{d}y \right| \leq C  \frac{\chi( \omega_\infty \Theta(t) ) \Theta(t)^2}{(t-\tau)^{1/2}} \int_\R e^{-\frac{|x-y|^2}{\kappa (t-\tau)}} (1+|y|)^2\omega(y)\md y,
\eqs
such for all $0<t<1$ we have
\bqs
\left| \int_0^t \int_\R \bG(t-\tau,x,y) \mathcal{N}(q_*(y),\omega(y)p(\tau,y))p(\tau,y)\mathrm{d}y\mathrm{d}\tau\right| \leq 2 C \chi( \omega_\infty \Theta(t) ) \Theta(t)^2 \int_\R (1+|y|)^2\omega(y)\md y,
\eqs
where we used that $\int_0^1(t-\tau)^{-1/2}\md \tau = 2$. As a consequence, one can find a positive constant (still denoted $C$) such that for all $0\leq t<1$ one has
\bqq
\Theta(t)\leq C ~  \| p_0\|_{L^\infty} + C \left( \int_\R (1+|y|)^2\omega(y)\md y\right) \chi( \omega_\infty \Theta(t) ) \Theta(t)^2. \label{eqzetasmall}
\eqq
Now for $t\geq 1$, using the estimate (\ref{eqGest}), we get
\bqs
\left| \int_\R \bG(t,x,y)p_0(y)\mathrm{d}y \right| \leq C \frac{1+|x|}{(1+t)^{3/2}}\int_\R(1+|y|)|p_0(y)|\md y,
\eqs
together with
\begin{align*}
\left| \int_\R \bG(t-\tau,x,y) \mathcal{N}(q_*(y),\omega(y)p(\tau,y))p(\tau,y)\mathrm{d}y \right| &\leq C  \frac{\chi( \omega_\infty \Theta(t) ) \Theta(t)^2(1+|x|)}{(1+t-\tau)^{3/2}(1+\tau)^3} \int_\R (1+|y|)^3\omega(y)\md y.
\end{align*}
Using the fact that \cite{xin92}
\bqs
\int_0^t \frac{1}{(1+t-\tau)^{3/2}(1+\tau)^3}\md \tau \leq \frac{\tilde{C}}{(1+t)^{3/2}},
\eqs
we obtain
\bqs
\left| \int_0^t \int_\R \bG(t-\tau,x,y) \mathcal{N}(q_*(y),\omega(y)p(\tau,y))p(\tau,y)\mathrm{d}y\mathrm{d}\tau\right| \leq \widehat{C} \frac{\chi( \omega_\infty \Theta(t) ) \Theta(t)^2(1+|x|)}{(1+t)^{3/2}}\int_\R (1+|y|)^3\omega(y)\md y.
\eqs
As a consequence, for all $t\geq 1$ we have
\bqq
\Theta(t)\leq C \int_\R(1+|y|)|p_0(y)|\md y+ \widehat{C} \left( \int_\R (1+|y|)^3\omega(y)\md y\right) \chi( \omega_\infty \Theta(t) ) \Theta(t)^2. \label{eqzeta}
\eqq
Finally, combing \eqref{eqzetasmall} and \eqref{eqzeta}, there exist constants $C_1>0$ and $C_2>0$ such that for all $t\in[0,T_*)$
\bqq
\Theta(t)\leq C_1 \Omega + C_2  \chi( \omega_\infty \Theta(t) ) \Theta(t)^2, \label{eqTheta}
\eqq
where we have set 
\[ \Omega:=\| p_0 \|_{L^\infty}+\int_{\mathbb{R}} (1+|y|)|p_0(y)|\mathrm{d}y. \]
So if we assume that the initial perturbation $p_0$ is small enough so that
\bqs
2C_1 \Omega <1, \quad \text{ and } 4C_1C_2 \chi( \omega_\infty 1 ) \Omega <1,
\eqs
then we claim that $\Theta(t)\leq 2C_1 \Omega < 1$ for all $t\in[0,T_*)$. Indeed, by eventually taking $C_1$ larger, we can always assume that $\Theta(0)=\sup_{x\in\mathbb{R}} \frac{|p_0(x)|}{1+|x|} \leq \|p_0\|_{L^\infty}<\Omega <2C_1 \Omega$ such that the continuity of $\Theta(t)$ implies that for small time we have $\Theta(t)< 2C_1 \Omega$. Suppose that there exists some $T>0$ where $\Theta(T)=2C_1 \Omega$ for the first time, then we have from \eqref{eqTheta}
\bqs
\Theta(T)\leq C_1\Omega\left(1+4C_1C_2 \chi( \omega_\infty 1 ) \Omega \right)<2C_1 \Omega,
\eqs
which is a contradiction and the claim is proved. This implies the maximal time of existence is $T_*=+\infty$ and the solution $p$ of \eqref{evolp} satisfies 
\bqs
\underset{t\geq 0}{\sup}~\underset{x\in\R}{\sup} ~ (1+t)^{3/2} \frac{|p(t,x)|}{1+|x|} < 2C_1 \Omega .
\eqs
This concludes the proof.

\section*{Acknowledgments} 
The authors thank the referees for comments that improved the paper.  
The authors would also like to thank the CIMI Excellence Laboratory, Toulouse, France, for
inviting MH as a Scientific Expert during the month of October 2017.  GF received support from the ANR project NONLOCAL ANR-14-CE25-0013.  MH received partial support from the National Science Foundation through grant NSF-DMS-1516155.  This work was partially supported by ANR-11-LABX-0040-CIMI within the program ANR-11-IDEX-0002-02.

\bibliographystyle{abbrv}
\bibliography{KPPBib}

\def\cprime{$'$}
\begin{thebibliography}{10}

\bibitem{alexander90}
J.~Alexander, R.~Gardner, and C.~Jones.
\newblock A topological invariant arising in the stability analysis of
  travelling waves.
\newblock {\em J. Reine Angew. Math.}, 410:167--212, 1990.

\bibitem{AW78}
D.~G. Aronson and H.~F. Weinberger.
\newblock Multidimensional nonlinear diffusion arising in population genetics.
\newblock {\em Advances in Mathematics}, 30(1):33--76, 1978.

\bibitem{beck14}
M.~Beck, T.~T. Nguyen, B.~Sandstede, and K.~Zumbrun.
\newblock Nonlinear stability of source defects in the complex
  {G}inzburg-{L}andau equation.
\newblock {\em Nonlinearity}, 27(4):739--786, 2014.

\bibitem{bricmont92}
J.~Bricmont and A.~Kupiainen.
\newblock Renormalization group and the {G}inzburg-{L}andau equation.
\newblock {\em Comm. Math. Phys.}, 150(1):193--208, 1992.

\bibitem{eckmann94}
J.-P. Eckmann and C.~E. Wayne.
\newblock The nonlinear stability of front solutions for parabolic partial
  differential equations.
\newblock {\em Comm. Math. Phys.}, 161(2):323--334, 1994.

\bibitem{focant98}
S.~Focant and T.~Gallay.
\newblock Existence and stability of propagating fronts for an autocatalytic
  reaction-diffusion system.
\newblock {\em Phys. D}, 120(3-4):346--368, 1998.

\bibitem{gallay94}
T.~Gallay.
\newblock Local stability of critical fronts in nonlinear parabolic partial
  differential equations.
\newblock {\em Nonlinearity}, 7(3):741--764, 1994.

\bibitem{gardner98}
R.~A. Gardner and K.~Zumbrun.
\newblock The gap lemma and geometric criteria for instability of viscous shock
  profiles.
\newblock {\em Comm. Pure Appl. Math.}, 51(7):797--855, 1998.

\bibitem{howard02}
P.~Howard.
\newblock Pointwise estimates and stability for degenerate viscous shock waves.
\newblock {\em J. Reine Angew. Math.}, 545:19--65, 2002.

\bibitem{johnson11}
M.~A. Johnson and K.~Zumbrun.
\newblock Nonlinear stability of spatially-periodic traveling-wave solutions of
  systems of reaction-diffusion equations.
\newblock {\em Ann. Inst. H. Poincar\'e Anal. Non Lin\'eaire}, 28(4):471--483,
  2011.

\bibitem{kapitula98}
T.~Kapitula and B.~Sandstede.
\newblock Stability of bright solitary-wave solutions to perturbed nonlinear
  {S}chr\"odinger equations.
\newblock {\em Phys. D}, 124(1-3):58--103, 1998.

\bibitem{kirchgassner92}
K.~Kirchgässner.
\newblock On the nonlinear dynamics of travelling fronts.
\newblock {\em Journal of Differential Equations}, 96(2):256 -- 278, 1992.

\bibitem{li16}
Y.~Li.
\newblock Point-wise stability of reaction diffusion fronts.
\newblock {\em arXiv preprint arXiv:1602.08176}, 2016.

\bibitem{raugel98}
G.~Raugel and K.~Kirchg\"assner.
\newblock Stability of fronts for a {KPP}-system. {II}. {T}he critical case.
\newblock {\em J. Differential Equations}, 146(2):399--456, 1998.

\bibitem{RW01}
V.~Rottsch{\"a}fer and C.~Wayne.
\newblock Existence and stability of traveling fronts in the extended
  {F}isher--{K}olmogorov equation.
\newblock {\em Journal of Differential Equations}, 176(2):532--560, 2001.

\bibitem{sandstede04}
B.~Sandstede and A.~Scheel.
\newblock Evans function and blow-up methods in critical eigenvalue problems.
\newblock {\em Discrete Contin. Dyn. Syst.}, 10(4):941--964, 2004.

\bibitem{sattinger}
D.~H. Sattinger.
\newblock On the stability of waves of nonlinear parabolic systems.
\newblock {\em Advances in Math.}, 22(3):312--355, 1976.

\bibitem{xin92}
J.~X. Xin.
\newblock Multidimensional stability of traveling waves in a bistable
  reaction--diffusion equation, i.
\newblock {\em Communications in partial differential equations},
  17(11-12):1889--1899, 1992.

\bibitem{zumbrun98}
K.~Zumbrun and P.~Howard.
\newblock Pointwise semigroup methods and stability of viscous shock waves.
\newblock {\em Indiana Univ. Math. J.}, 47(3):741--871, 1998.

\end{thebibliography}

\end{document}